\documentstyle[psfig, 11pt]{article}
\title{Monotone drawings of planar graphs\footnote{This is a revised version of \cite{PT04}. See the Remark at the end of Section 2.}}

\author{{\sl J\'anos Pach}\thanks{Supported by NSF grant
CR-00-98246, PSC-CUNY Research Award 63382-0032 and OTKA-T-032452.} \\
City College, CUNY and \\ Hungarian Academy of Sciences\\
\and
{\sl G\'eza T\'oth}\thanks{Supported by NSF grant OTKA-T-038397,
OTKA-T-032452.}\\ R\'enyi Institute of the\\ Hungarian
Academy of Sciences\\}

\date{}
\begin{document}
\maketitle

\begin{abstract}
Let $G$ be a graph drawn in the plane so that its edges are
represented by $x$-monotone curves, any pair of which cross
an even number of times. We show that $G$ can be redrawn in
such a way that the $x$-coordinates of the vertices remain
unchanged and the edges become non-crossing straight-line
segments.
\end{abstract}

\section{Introduction}

A {\em drawing} ${\cal D}(G)$ of a graph $G$ is a representation
of the vertices and the edges of $G$ by points and by possibly
crossing simple Jordan arcs connecting the corresponding
point pairs, resp. When it does not lead to confusion, we make
no notational or terminological distinction between the vertices
(resp. edges) of the underlying abstract graph and the points
(resp. arcs) representing them. Throughout this paper, we assume
that in a drawing
\begin{enumerate}
\item  no edge passes through any vertex other than its
endpoints;
\item  any two edges cross only a finite number of times;
\item  no three edges cross at the same point;
\item  if two edges of a drawing share an interior
point $p$ then they properly cross at $p$, i.e., one arc
passes from one side of the other arc to the other side;
\item  no two vertices have the same $x$-coordinate.
\end{enumerate}

A drawing is called {\em $x$-monotone} if every vertical
line intersects every edge in at most one point.
We call a drawing {\em even} if any two edges
cross an even number of times.

Hanani (Chojnacki) \cite{Ch34} (see also \cite{T70})
proved the remarkable theorem that if a graph $G$ permits an
even drawing, then it is {\em planar}, i.e., it can be redrawn
without any crossing. On the other hand, by F\'ary's theorem
\cite{F48}, \cite{W36}, every planar graph has a straight-line
drawing. We can combine these two facts by saying that
every even drawing can be {\em ``stretched''}.

The aim of this note is to show that if we restrict our attention to
{\em $x$-monotone} drawings, then every even drawing can
be stretched without changing the $x$-coordinates of the vertices.

Consider an $x$-monotone drawing ${\cal D}(G)$ of a graph $G$. If
the vertical ray starting at $v\in V(G)$ and pointing upward
(resp. downward) crosses an edge $e\in E(G)$, then $v$ is said to
be {\em below} (resp. {\em above}) $e$. Two drawings of the same
graph are called {\em equivalent}, if in a small neighborhood
of each vertex $v\in V(G)$, the above-below relationships between the edges adjacent to $v$ are the same.

In the next two sections we establish the following two results.

\medskip

\noindent {\bf Theorem 1.} {\em For any $x$-monotone even
drawing of a connected graph, there is an equivalent $x$-monotone
drawing in which no two edges cross each other and the
$x$-coordinates of the corresponding vertices are the same.}

\medskip

\noindent {\bf Theorem 2.} {\em For any non-crossing
$x$-monotone drawing of a graph $G$, there is an equivalent
non-crossing straight-line drawing, in which the
$x$-coordinates of the corresponding vertices are the same.}

\medskip

Two edges are called {\em adjacent} if they share an endpoint.
It is an interesting open problem to decide whether Theorem 1
remains true under the weaker assumption that any two
{\em non-adjacent} edges cross an even number of times. Hanani's
theorem mentioned above is valid in this stronger form. It
was suggested by Tutte ``that crossings of adjacent edges
are trivial, and easily got rid of.'' We have been unable to
verify this view.

\section{Proof of Theorem 1}

We can assume that $G$ is connected.
We follow the approach of Cairns and Nikolayevsky \cite{CN00}.
Consider an $x$-monotone drawing $\cal D$ of a graph on the
$xy$-plane, in which any two edges cross an even number of times.
Let $u$ and $v$ denote the leftmost and rightmost vertex,
respectively. We can assume without loss of generality that
$u=(-1, 0)$ and $v=(1, 0)$. Introduce two additional vertices,
$w=(0,1)$ and $z=(0,-1),$ each connected to $u$ and $v$ by arcs
of length $\pi/2$ along the unit circle $C$ centered at the
origin, and suppose that every other edge of the drawing lies in
the interior of $C$. Denote by $G$ the underlying abstract graph,
including the new vertices $w$ and $z$.


For each crossing point $p$, attach a {\em handle} (or bridge)
to the plane in a very small neighborhood $N(p)$ of $p$, with
radius $\varepsilon>0$. Assume that (1) these neighborhoods
are pairwise disjoint, (2) $N(p)$ is disjoint from every other
edge that does not pass through $p$, and that (3) every vertical
line intersects every handle only at most once.  For every $p$,
take the portion belonging to $N(p)$ of one of the edges that
participate in the crossing at $p$, and lift it to the handle
without changing the $x$- and $y$-coordinates of its points.
The resulting drawing ${\cal D}_0$ is a crossing-free
embedding of $G$ on a surface $S_0$ of possibly higher genus.

Let $S_1$ be a very small closed neighborhood of the drawing
${\cal D}_0$ on the surface $S_0$, with positive radius
$\varepsilon'<\varepsilon$. Note that $S_1$ is a compact,
connected surface, whose boundary consists of a finite number
of closed curves. Attaching a disk to each of these closed
curves, we obtain a surface $S_2$ with no boundary. According
to Cairns and Nikolayevsky \cite{CN00}, $S_2$ must be a
$2$-dimensional {\em sphere}. To verify this claim, consider
two closed curves, $\alpha_2$ and $\beta_2$, on $S_2$. They
can be deformed into closed walks, $\alpha_1$ and $\beta_1$,
respectively, along the edges of ${\cal D}_0$. The projection
of these two walks into the $(x,y)$-plane are closed walks,
$\alpha$ and $\beta$ in $\cal D$, that must cross each other
an even number of times. Every crossing between $\alpha$ and
$\beta$ occurs either at a vertex of $\cal D$ or between two
of its edges. By the assumptions, any two edges in $\cal D$ cross
an even number of times. (The same assertion is trivially true in
${\cal D}_0\subset S_2$, because there no two edges cross.) Using
the fact that in ${\cal D}_0\subset S_2$ the cyclic order of the
edges incident to a vertex is the same as the cyclic order of the
corresponding edges in $\cal D$, we can conclude that $\alpha_1$
and $\beta_1$ cross an even number of times, and the same is true
for $\alpha_2$ and $\beta_2$. Thus, $S_2$ is a surface with no
boundary, in which any two closed curves cross an even number of
times. This implies that $S_2$ is a sphere. Consequently,
${\cal D}_0$, a crossing-free drawing of $G$ on $S_2$, corresponds
to a plane drawing.

For any point $q$, 
let $x(q)$ denote the $x$-coordinate of $q$. As before, every boundary curve of $S_1$ corresponds to a cycle of $G$. Since in the
original drawing, the cycle $vwuz$ encloses all other edges and
vertices of $G$, one of the boundary curves of $S_1$, say
$\gamma$, corresponds to the cycle $vwuz$. Let $D_{\gamma}$ be the disk attached to $\gamma$. Since $S_2$ is homeomorphic to a sphere, $S_2\setminus int(D_{\gamma})$ is homeomorphic to a closed disk $D$, whose boundary corresponds to the cycle  $vwuz$.
We will define a function $f$ on the points $p\in D$ such that $f(p)$ can be regarded as the ``$x$-coordinate of $p$." Using this function, the drawing ${\cal D}_0$ can also be regarded as an {\em
$x$-monotone} plane drawing of $G$, in which the $x$-coordinates
of the vertices are the same as the $x$-coordinates of the
corresponding vertices in $\cal D$.

First, we prove Theorem 1 for cycles.

\medskip

\noindent{\bf Lemma 2.1.}
{\em For any $x$-monotone even drawing of a cycle, there is an equivalent non-crossing straight-line drawing, in which the $x$-coordinates of the corresponding vertices are the same.}


\medskip

\noindent{\bf Proof.} Suppose that $C=v_1v_2, \cdots v_i$ is a
cycle with an $x$-monotone even drawing. For $i=3,4$, the lemma can be easily verified. Let $i>3$, and suppose that we have already proved the assertion for every integer smaller than $i$.
Let $x_1, x_2, \ldots , x_i$ denote the $x$-coordinates of
$v_1, v_2, \ldots , v_i$, respectively. Choose an index $j$ for which $|x_{j+1}-x_j|$ is {\em minimum}, where the indices are taken modulo $i$. Suppose without loss of generality that $x_j<x_{j+1}$. If we have $x_{j+1}<x_{j+2}$ (or $x_{j-1}>x_j$),
then delete $x_{j+1}$ (resp., $x_j$), apply the lemma to the
remaining sequence, and insert an extra vertex $v_{j+1}$ (resp. $v_j$) whose $x$-coordinate is $x_{j+1}$ (resp., $x_j$) in the corresponding side of the resulting polygon.  Otherwise, by the minimality assumption, we have $x_{j+2}<x_j<$, $x_{j+1}<x_{j-1}$.
In this case, apply the
lemma to the sequence obtained by the deletion of $v_j$ and
$v_{j+1}$, and notice that the $v_{j-1}v_{j+2}$
side of the resulting polygon,
whose endpoints have $x$-coordinates $x_{j-1}$ and $x_{j+2}$,
can be replaced by three edges meeting the requirements,
running very close to it.  $\Box$

\medskip

Consider the drawing ${\cal D}_0$ of $G$ on $S_1$. For each point $p$ on the
edges of $G$, let $f(p)=x(p)$, where $x(p)$ denotes the $x$-coordinate of $p$
in the original drawing 
$\cal D$.

Let $\kappa$ be a boundary curve of $S_1$, distinct from $\gamma$, the
boundary curve corresponding to the cycle $vwuz$.
Let $C_{\kappa}$ be the cycle of $G$ that corresponds to $\kappa$, as it was
drawn in the 
original drawing $\cal D$.
Apply Lemma 2.1 to $C_{\kappa}$, and denote the resulting drawing by $C'_{\kappa}$. Let $D_{\kappa}$
be the closed polygonal region (topological disk) bounded by $C'_{\kappa}$. For any point $p\in D_{\kappa}$ let $f(p)=x(p)$.
The points of $\kappa$ (a boundary curve of $S_1$) and the points of
$C'_{\kappa}$ (the boundary of $D_{\kappa}$) 
are both in one-to-one
correspondence with the points of $C_{\kappa}$.
Attach $D_{\kappa}$ to $\kappa$ so that the points attached to each other correspond to the same point of $C_{\kappa}$.
Repeating the same procedure for each boundary curve of $S_1$, different from
$\gamma$, we obtain a crossing-free drawing of $G$ on $D$, together with a
continuous function $f(q)$ defined on $D$, which coincides with $x(q)$ for
every point $q$ that lies on an edge or on a vertex of $G$.  By our
construction, we have $f(u)=-1$, $f(v)=1$, and $-1<f(q)<1$ for each $q\in D$, 
$q\ne u, v$.

In order to justify the claim that $f(q)$ can be regarded as the
$x$-coordinate of $q$ in the new drawing, we have to show that for any fixed
$x$, $-1<x<1$, the set $L(x)=\{\  q\ |\  f(q)=x\ \}$ is a simple curve connecting
two boundary points of $D$. Clearly, there is exactly one point $q_1$
(resp. $q_2$) on the path $uzv$ (resp. $uwv$) with $x(q_1)=f(q_1)=x$ 
(resp. $x(q_2)=f(q_2)=x$.
If $L(x)$ is not a level curve connecting $q_1$ and $q_2$,
then it must contain a loop (a simple closed subcurve). In the interior of
such a loop, $f$ must have a local maximum or minimum, say, at a point
$r$. Thus, it is enough to show that no such $r$ exists. If $r$ lies in the
interior of a disk $D_{\kappa}$, then it cannot be locally extreme, because in
such a region $f$ is defined as the $x$-coordinate of the points in a planar
embedding of $D_{\kappa}$. If $r$ lies in the interior of an edge, then it
cannot be locally extreme either, since restricted to edges, $f$ is a strictly
monotone function. We are left with the case when 
$r$
is a vertex of $G$. If there is at least one edge incident to $r$ on both
sides of $r$, then we can argue in the same way as in the last 
case.

The only remaining case is when $r$ is a vertex and all edges incident to $r$
are on one side of $r$. To deal with this case, we 
need some preparation.

Let $C=v_1v_2\cdots v_i$ be a cycle (closed curve) in the plane, passing
through the points $v_i$ in this order. Orient it arbitrarily. Given a point
$p$ not on $C$, its {\em winding number} $w(p)$ is the number of times $C$
travels counterclockwise around $p$. The {\em interior} $I(C)$ and the {\em
  exterior} $E(C)$ of $C$ are defined as the set of all points in the plane
with odd winding number and the set of all points with even winding number,
respectively. If we reverse the orientation of $C$, its interior and the
exterior remain unchanged. Apart from a bounded region, all points of the
plane belong to the exterior of $C$.

Let $v_j$ be one of the vertices of $C$. The edges (arcs) $v_jv_{j-1}$ and $v_jv_{j+1}$ divide a small neighborhood of $v_j$
into two parts; one of them belongs to $I(C)$, the other to $E(C)$.
Listing the arcs and regions in the counter-clockwise order around $v_j$, there are two possibilities:
$v_jv_{j-1}, I(C), v_jv_{j+1}, E(C)$, or
$v_jv_{j-1}, E(C), v_jv_{j+1}, I(C)$. In the first case,
$v_j$ is said to be of {\em type 1}, in the second case
it is said to be of {\em type 2}.

\medskip

\noindent{\bf Lemma 2.2.} {\em Let $C$ and $C'$ be two equivalent
$x$-monotone even drawings of a cycle $v_1v_2\cdots v_i$, in which
the $x$-coordinates of the corresponding vertices are the same.
Then the type of each vertex is the same in both drawings.}

\medskip

\noindent{\bf Proof.} Suppose that $v_1$ is the leftmost vertex of
$C$. Then in both drawings, both $v_1v_i$ and $v_1v_2$ lie to the
right of $v_1$. Assume without loss of generality that in $C$, in
a small neighborhood of $v_1$, the arc $v_1v_i$ is {\em below}
$v_1v_2$. Since $v_1$ is the leftmost vertex, $I(C)$ must lie to
the right of $v_1$. Thus, in $C$, vertex $v_1$ is of type 1. It
follows from the equivalence of the two drawings that in $C'$, in
a small neighborhood of $v_1$, the arc $v_1v_i$ lies below
$v_1v_2$ and $I(C')$ is to the right of $v_1$. Hence, in $C'$ the
vertex $v_1$ is also of type 1. In particular, in both drawings,
in a small neighborhood of $v_1$, region $I(C)$, resp. $I(C')$,
must lie {\em below} $v_1v_2$. Moving from $v_1$ to $v_2$ along
the edge $v_1v_2$, we encounter an even number of crossings.
Therefore, in both drawings, in a small neighborhood of $v_2$, the
region $I(C)$, resp. $I(C')$, also lies {\em below} $v_1v_2$.
This, in turn, implies that the type of $v_2$ is also the same in
both drawings. In the same way, we can prove by induction that the
types of $v_3, \ldots , v_i$ are the same in both drawings. $\Box$

\medskip

Return to the proof of Theorem 1. We were left with the case,
where $r$ is a vertex of $G$ and all edges incident to $r$ are on
the same side of $r$, say, to the {\em left} of it. We will show
that the function $f$ cannot attain a local extremum at $r$.
Obviously, it cannot attain a local {\em minimum}.

Consider a small neighborhood of $r$. Let $e_1, e_2, \ldots , e_i$
denote the edges incident to $r$, listed in counter-clockwise
order around $r$. For any $j$, $1\le j\le i$, let $\kappa_j$
denote the uniquely determined boundary curve of $S_1$, in which
the arcs corresponding to $e_j$ and $e_{j+1}$ are consecutive.
(The indices are taken modulo $i$.) Let $C_j$ denote the cycle in
$G$ which corresponds to $\kappa_j$ in the original drawing $\cal
D$. Using our notation, we have $C_j=C_{\kappa_j}$.

We claim that in a small neighborhood of $r$, $\kappa_j$ is in the
{\em interior} of $C_j$. Notice that this claim is true if and
only if $\kappa_j$ is in the {\em interior} of $C_j$ in a small
neighborhood of any other vertex of $C_j$. (This follows from the
fact that $\cal D$ is an even drawing and $\kappa_j$ is a boundary
curve of $S_1$.) Since $r\ne u, w, v, z$, we have
$\kappa\ne\gamma$, that is, $C_j\ne uwvz$. Therefore, in $\cal D$,
at least one of the vertices $u, w, v, z$ is in the {\em exterior}
of $C_j$. Take such a vertex and a shortest path connecting it to
a vertex $p_1$ of $C_j$. Let $p_2$ be the previous vertex along
this path. Clearly, $p_2$ belongs to the exterior of $C_j$,
because any two edges cross an even number of times. In a small
neighborhood of $p_1$, $\kappa_j$ lies between two consecutive
edges incident to $p_1$, so the edge $p_1p_2$ lies on the side of
$C_j$ opposite to $\kappa_j$. Since $p_2$ belongs to the exterior
of $C_j$, and $p_1p_2$ crosses $C_j$ an even number of times, in a
small neighborhood of $p_1$, the edge $p_1p_2$ is in the exterior
and $\kappa$ in the {\em interior} of $C_j$.

Consider now $C'_j$, the crossing-free drawing of $C_j$, meeting
the requirements of Lemma 2.1. We glued $D_j$, the {\em interior}
of $C'_j$, to $\kappa_j$, and repeated this procedure for every
$j$. Consider now the index $j$, for which the interior of $C_j$
contains a short horizontal segment whose left endpoint is $r$.
Starting at $r$ and moving along this segment to the right, the
$x$-coordinates of the points increase. Applying Lemma 2.2 to
$C_j$ and $C'_j$, we can conclude that starting at $r$, within
$D_j$ we can also move to the right. Therefore, along such a path
$f$ increases. This implies that $r$ cannot attain a local {\em
maximum} at $r$.

Summarizing: ${\cal D}_0$ is a crossing-free drawing of $G$ in a
disc $D$, and $f$ is a function defined on $D$. Along the vertices
and edges of $G$, $f$ was defined to be equal to the
$x$-coordinate of the corresponding point in the original drawing
${\cal D}$. Each level curve of $f$ is a simple curve connecting a
pair of boundary points of $D$. Therefore, the level curves can be
consistently parameterized so that the new parameter can be
regarded as the $y$-coordinate, and the function $f$ as the
$x$-coordinate of the points. The resulting drawing satisfies the
requirements of Theorem 1.
 $\Box$.

\medskip

\noindent {\bf Remark.} We are grateful to M. Pelsmajer and M. Schaefer, who
pointed out a mistake in the published version of the above proof. Originally,
we defined two drawings to be equivalent if the above-below relationship
between vertices and edges are the same. However, one can guarantee only the
weaker property that in the new drawing the above-below relationship is
preserved in small neighborhoods of the vertices. In the present version, two
$x$-monotone drawings are defined to be equivalent if they satisfy this
condition.

\section{Proof of Theorem 2}

Let ${\cal D}={\cal D}(G)$ be a non-crossing $x$-monotone
drawing of a graph $G$. First, we show that it is
sufficient to prove Theorem 2 for {\em triangulated} graphs.
Deleting all vertices (points) and edges (arcs) of
${\cal D}$ from the plane, the plane falls into connected
components, called {\em faces}. The $x$-coordinate of any
vertex $v$ will be denoted by $x(v)$.

\medskip

\noindent {\bf Lemma 3.1.}
{\em By the addition of further edges
and an extra vertex, if necessary, every
non-crossing $x$-monotone drawing $\cal D$ can be extended to
a non-crossing $x$-monotone triangulation.}

\medskip

\noindent {\bf Proof.} Consider a face $F$, and assume that
it has more than $3$ vertices. It is sufficient to show that
one can always add an $x$-monotone edge between two
non-adjacent vertices of $F$, which does not cross any
previously drawn edges.

For the sake of simplicity, we outline the argument
only for the case when $F$ is a bounded face.
The proof in the other case is very similar, the only
difference is that we may also have to add an extra vertex.

\vskip 0.6cm

\centerline{\psfig{figure=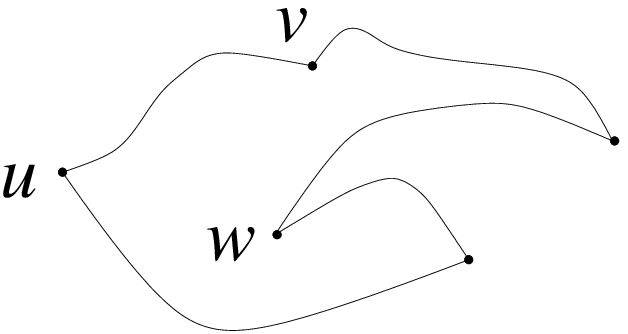}}

\bigskip

\smallskip

\centerline{{\bf Figure 1.} {\em The vertex $w$ is extreme,
$u$ and $v$ are not.}}

\medskip

A vertex $w$ of $F$ is called {\em extreme} if it is not
the left endpoint of any edge or not the right endpoint
of any edge in $\cal D$, and a small neighborhood of $w$
on the vertical line through $w$ belongs to $F$. In
particular, if the boundary of $F$ is not connected, the
leftmost (and the rightmost) vertex of each component of
the boundary other than the exterior component, is extreme.
See Fig. 1.

Suppose first that $F$ has an extreme vertex $w$. We may
assume, by symmetry, that $w$ is not the right endpoint
of any edge in $\cal D$. Starting at $w$, draw a horizontal
ray in the direction of the negative $x$-axis. Let $p$ be
the first intersection point of this ray with the boundary
of $F$. If $p$ is a vertex, then the segment $wp$ can be
added to $\cal D$. Otherwise, one can add an $x$-monotone
edge joining $w$ to the left endpoint of the edge that $p$
belongs to.

Suppose next that none of the vertices of $F$ are extreme.
In this case, the boundary of $F$ is connected and any two
vertices of $F$ can be joined by an $x$-monotone curve
inside $F$. However, an edge can be added to $\cal D$ only
if the corresponding two vertices do not induce an edge in
the exterior of $F$. Clearly, letting $v_1$, $v_2$, $v_3$,
and $v_4$ denote four consecutive vertices of $F$, at least
one of the pairs $(v_1,v_3)$ and $(v_2,v_4)$ has this
property. $\Box$

\medskip

Now we turn to the proof of Theorem 2. The proof is by
induction on the number of vertices. If $G$ has at most $4$
vertices, the assertion is trivial. Suppose that $G$ has
$n>4$ vertices and that we have already established the
theorem for graphs having fewer than $n$ vertices. By Lemma
3.1, we can assume without loss of generality that the
original $x$-monotone drawing $\cal D$ of $G$ is
triangulated.

\medskip

\noindent {\sc Case 1.} There is a triangle $T=v_1v_2v_3$ in
$\cal D$, which is not a face.

Then there is at least one vertex of $\cal D$ in the
interior and at least one vertex in the exterior of $T$.
Consequently, the drawings ${\cal D}_{\mbox {\scriptsize in}}$
and ${\cal D}_{\mbox {\scriptsize out}}$ defined as the part of
${\cal D}$ induced by $v_1, v_2$, $v_3$, and all vertices
{\em inside} $T$ and {\em outside} $T$, resp., have fewer than
$n$ vertices. By the induction hypothesis, there exist
straight-line drawings ${\cal D}'_{\mbox {\scriptsize in}}$
and ${\cal D}'_{\mbox {\scriptsize out}}$, equivalent to
${\cal D}_{\mbox {\scriptsize in}}$ and
${\cal D}_{\mbox {\scriptsize out}}$, resp., in which all
vertices have the same $x$-coordinates as in the original
drawing. Notice that there is an affine transformation $A$
of the plane, of the form
$$A(x,y)=(x,ax+by+c),$$
which takes the triangle induced by $v_1, v_2$, $v_3$ in
${\cal D}_{\mbox {\scriptsize in}}$ into the triangle induced
by $v_1, v_2$, $v_3$ in ${\cal D}_{\mbox {\scriptsize out}}$.
Since the image of a drawing under any affine transformation
is equivalent to the original drawing, we conclude that
$A\left({\cal D}'_{\mbox {\scriptsize in}}\right)\cup
{\cal D}'_{\mbox {\scriptsize out}}$ meets the requirements.

\medskip

In the sequel, we can assume that $\cal D$ has no triangle that
is not a face. Fix a vertex $v$ of $\cal D$ with minimum
degree. Since every triangulation on $n>4$ vertices has
$3n-6$ edges, the degree of $v$ is $3, 4,$ or $5$. If the
degree of $v$ is $3$, the neighbors of $v$ induce a triangle
in $\cal D$, which is not a face, contradicting our
assumption.

There are two more cases to consider.

\bigskip

\noindent {\sc Case 2.} The degree of $v$ is 4.

Let $v_1, v_2, v_3, v_4$ denote the neighbors of $v$, in
clockwise order. There are three substantially different
subcases, up to symmetry. See Fig. 2.

\vskip 0.6cm

\centerline{\psfig{figure=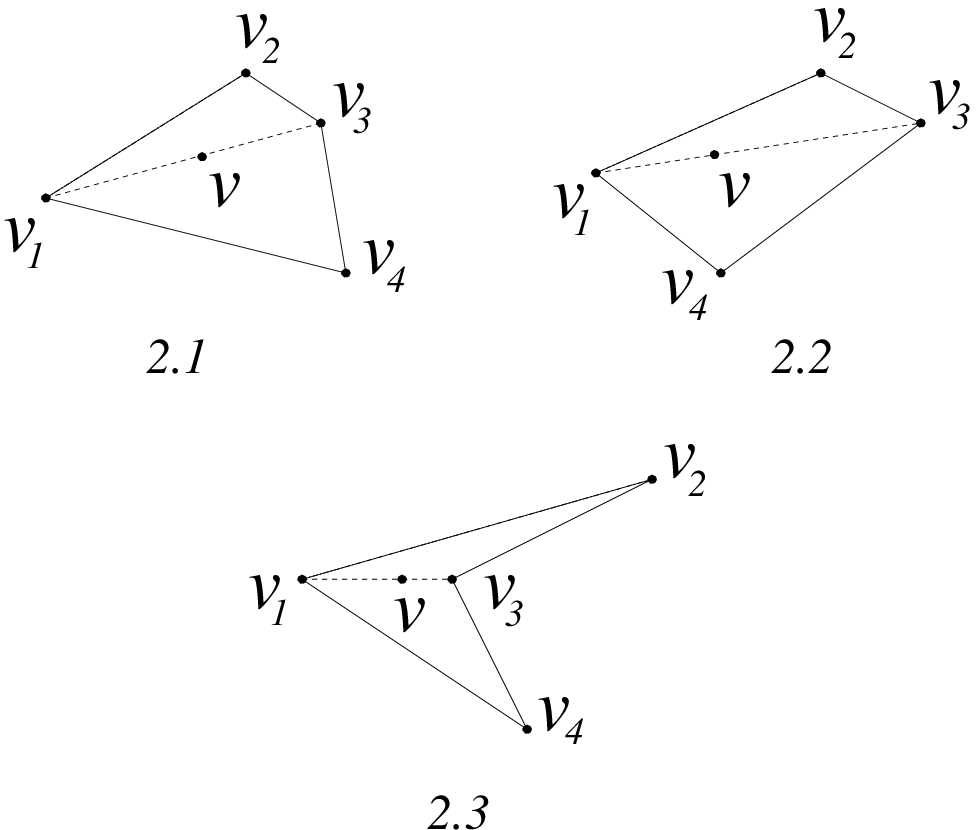}}

\bigskip

\smallskip

\centerline{{\bf Figure 2.} {\sc Case 2.}}

\medskip

{\sc Subcase 2.1:} $x(v_1)<x(v_2)<x(v_3)<x(v_4)$

Clearly, at least one of the inequalities $x(v)>x(v_2)$ and
$x(v)<x(v_3)$ is true. Suppose without loss of generality
that $x(v)<x(v_3)$. If $v_1$ and $v_3$ were connected by an
edge, then $vv_1v_3$ would be a triangle with $v_2$ and
$v_4$ in its interior and in its exterior, resp.,
contradicting our assumption. Remove $v$ from $\cal D$, and
add an $x$-monotone edge between $v_1$ and $v_3$, running
in the interior of the face that contains $v$. Applying the
induction hypothesis to the resulting drawing, we obtain
that it can be redrawn by straight-line edges, keeping the
$x$-coordinates fixed. Subdivide the segment $v_1v_3$ by
its (uniquely determined) point whose $x$-coordinate is
$x(v)$. In this drawing, $v$ can also be connected by
straight-line segments to $v_2$ and to $v_4$. Thus, we
obtain an equivalent drawing which meets the requirements.

\smallskip

{\sc Subcase 2.2:} $x(v_1)<x(v_2)<x(v_3)>x(v_4)>x(v_1)$

{\sc Subcase 2.3:} $x(v_1)<x(v_2)>x(v_3)<x(v_4)>x(v_1)$

In these two subcases, the above argument can be repeated {\em
verbatim}. In Subcase 2.3, to see that $x(v_1)<x(v)<x(v_3)$, we
have to use the fact that in $\cal D$ both $vv_2$ and $vv_4$ are
represented by $x$-monotone curves.

\medskip

\noindent {\sc Case 3.} The degree of $v$ is 5.

Let $v_1, v_2, v_3, v_4, v_5$ be the neighbors of $v$, in
clockwise order. There are four substantially different
cases, up to symmetry. See Fig. 3.

\smallskip

{\sc Subcase 3.1:} $x(v_1)<x(v_2)<x(v_3)<x(v_4)<x(v_5)$

{\sc Subcase 3.2:} $x(v_1)<x(v_2)<x(v_3)<x(v_4)>x(v_5)>x(v_1)$

{\sc Subcase 3.3:} $x(v_1)<x(v_2)<x(v_3)>x(v_4)<x(v_5)>x(v_1)$

{\sc Subcase 3.4:} $x(v_1)<x(v_2)>x(v_3)>x(v_4)<x(v_5)>x(v_1)$

In all of the above subcases, we can assume, by symmetry
or by $x$-monotonicity, that $x(v)<x(v_4)$.
Since $\cal D$ has no triangle which is
not a face, we obtain that $v_1v_3$, $v_1v_4$, and $v_2v_4$
cannot be edges. Delete from $\cal D$ the vertex $v$ together
with the five edges incident to $v$, and let ${\cal D}_0$
denote the resulting drawing.
Furthermore, let ${\cal D}_1$ (and ${\cal D}_2$) denote the
drawing obtained from ${\cal D}_0$ by adding two
non-crossing $x$-monotone diagonals, $v_1v_3$ and
$v_1v_4$ (resp. $v_2v_4$ and $v_1v_4$), which run in the
interior of the face containing $v$. By the induction
hypothesis, there exist straight-line drawings
${\cal D}'_1$ and ${\cal D}'_2$ equivalent to ${\cal D}_1$
and ${\cal D}_2$, resp., in which the $x$-coordinates of
the corresponding vertices are the same.

Apart from the edges $v_1v_3$, $v_1v_4,$ and $v_2v_4$,
${\cal D}'_1$ and ${\cal D}'_2$ are non-crossing straight-line
drawings equivalent to ${\cal D}_0$ such that the
$x$-coordinates of the corresponding vertices are the same.
Obviously, the convex combination of two such drawings is
another non-crossing straight-line drawing equivalent to
${\cal D}_0$. More precisely, for any $0\le \alpha \le 1$, let
${\cal D}'_{\alpha}$ be defined as
$${\cal D}'_{\alpha}=\alpha{\cal D}'_1+(1-\alpha ){\cal D}'_2.$$
That is, in ${\cal D}'_{\alpha}$, the $x$-coordinate of any
vertex $u\in V(G)-v$ is equal to $x(u)$, and its $y$-coordinate
is the combination of the corresponding $y$-coordinates in
${\cal D}'_1$ and ${\cal D}'_2$ with coefficients $\alpha$ and
$1-\alpha$, resp.

Observe that the only possible concave angle of the quadrilateral
$Q=v_1v_2v_3v_4$ in ${\cal D}'_1$ and ${\cal D}'_2$ is at $v_3$
and at $v_2$, resp. In ${\cal D}'_{\alpha}$, $Q$ has at most one
concave vertex. Since the shape of $Q$ changes continuously with
$\alpha$, we obtain that there is a value of $\alpha$ for which
$Q$ is a {\em convex} quadrilateral in ${\cal D}_{\alpha}$. Let
${\cal D}'$ be the straight-line drawing obtained from ${\cal
D}'_{\alpha}$ by adding $v$ at the unique point of the segment
$v_1v_4$, whose $x$-coordinate is $x(v)$, and connect it to
$v_1,\ldots,v_5$. Clearly, ${\cal D}'$ meets the requirements of
Theorem 2.

\vskip 1.6cm

\centerline{\psfig{figure=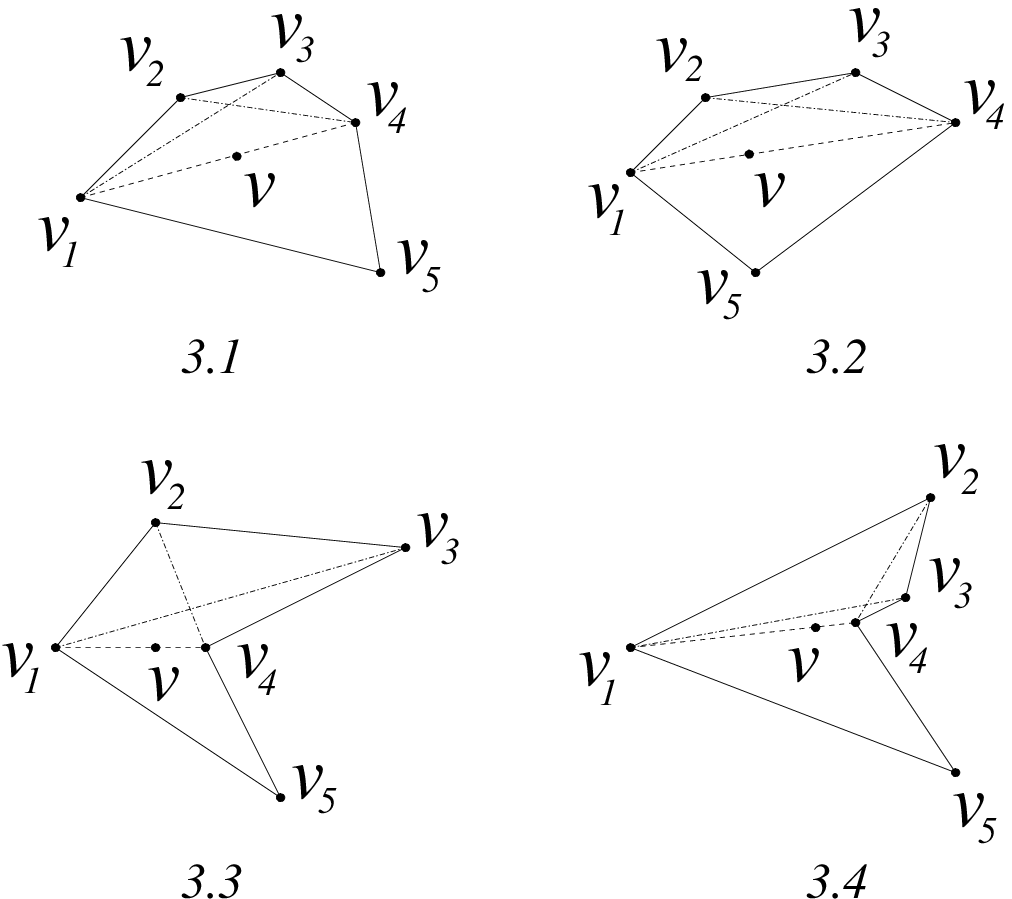}}

\bigskip

\smallskip

\centerline{{\bf Figure 3.} {\sc Case 3.}}

\bigskip

\noindent {\bf Remark:} We are grateful to Professor
P. Eades for calling our attention to his paper
\cite{EFL96}, sketching a somewhat more complicated proof
for a result essentially equivalent to our Theorem 2.

\end{document}